\documentclass[11pt]{article}
\usepackage[a4paper,margin=27mm]{geometry}
\usepackage[T1]{fontenc}
\usepackage{lmodern}
\usepackage{amsmath,amssymb,amsthm,mathtools,bm}
\usepackage{microtype}
\usepackage{hyperref}
\hypersetup{colorlinks=true,linkcolor=blue,citecolor=blue,urlcolor=blue,
  pdfauthor={Shuyuan Liang},
  pdftitle={Wasserstein bounds and Cramer-type moderate deviations for the multi-dimensional elephant random walk}}
\numberwithin{equation}{section}
\allowdisplaybreaks[2]
\newtheorem{theorem}{Theorem}[section]
\newtheorem{lemma}[theorem]{Lemma}
\newtheorem{proposition}[theorem]{Proposition}
\theoremstyle{definition}
\newtheorem{remark}[theorem]{Remark}
\newcommand{\E}{\mathbf E}
\newcommand{\Pp}{\mathbf P}
\newcommand{\R}{\mathbb R}
\newcommand{\F}{\mathcal F}
\newcommand{\Lw}{\mathcal L}
\newcommand{\op}{\mathrm{op}}
\newcommand{\tr}{\operatorname{tr}}
\newcommand{\Cov}{\operatorname{Cov}}
\newcommand{\ind}{\mathbf 1}
\newcommand{\Wp}{W_{p,2}}
\newcommand{\norm}[1]{\left\lVert #1\right\rVert}

\title{$p$-Wasserstein distance and Cram\'{e}r-type moderate deviations
for the multi-dimensional elephant random walk}
\author{Shuyuan Liang\thanks{Corresponding author. Email:
\href{mailto:liangshuyuan0501@foxmail.com}{\nolinkurl{liangshuyuan0501@foxmail.com}}}\\[0.5ex]
\small School of Mathematics and Statistics, Northeastern University at Qinhuangdao,\\
\small Qinhuangdao, {\rm066004}, Hebei, China.}
\date{}
\begin{document}
\maketitle

\begin{abstract}
In this paper, we establish $p$-Wasserstein bounds and Cram\'{e}r-type moderate
deviations for the multi-dimensional elephant random walk in the diffusive and
critical regimes. Using a Gaussian convolution estimate with random covariance
and the recursions for the conditional covariance, we obtain
$W_{p,2}(\Lw(\widehat T_n),N(0,I_d))\le C_{a,d}pL_n$, where $\widehat T_n$
is normalized by its exact covariance and
\[
L_n\asymp
\begin{cases}
n^{-1/2}\log n,&a\le1/3,\\
n^{-3(1/2-a)},&1/3<a<1/2,\\
(\log n)^{-3/2},&a=1/2.
\end{cases}
\]
Consequently, we obtain radial Cram\'{e}r-type moderate deviations for
$x=o(L_n^{-1/3})$. We also determine the covariance bias under the uncorrected
normalization and establish the corresponding Wasserstein bounds, lower bounds
and moderate deviation estimates. The dimension and the model parameters are fixed.
\end{abstract}

\begingroup\small
\noindent\textbf{Keywords.} Multi-dimensional elephant random walk; Martingale central limit theorem;
Wasserstein distance; Random covariance; Cram\'er-type moderate deviations

\noindent\textbf{MSC2020.} 60F05, 60F10, 60G42, 60G50
\par\endgroup

\section{Introduction}\label{sec:introduction}
\quad
The elephant random walk (ERW) is a random walk with complete memory of its past,
first introduced by Sch\"utz and Trimper~\cite{schutz2004elephants}.
Its asymptotic properties have attracted considerable attention, see \cite{BMERW, B21, CGSCLT, ASIP, FanHuMa2021, MEFWD}.
Moreover, Bercu and Laulin~\cite{BL} introduced the multi-dimensional elephant random walk (MERW) and established
multivariate limit theorems by a martingale approach.

Deviation estimates have also been studied for the ERW and its extensions.
Fan and Shao~\cite{FanShao2024} established Cram\'er-type moderate deviation
expansions for martingales, with applications to the ERW. For the ERW with
random step sizes, Fan \textit{et al.}~\cite{FanLiangMa2026} studied
Cram\'er-type moderate deviations and precise large deviations, while
Liang~\cite{Liang2026} obtained local precise large deviations.

The MERW can be described as follows. For $d\ge1$, set $(e_j)_{1\leq j\leq d}\in\mathbb Z^d$ and 
\[
e_j=(0,\ldots,0,\underset{\substack{\uparrow\\ j}}{1},0,\ldots,0)^{\top},\qquad 1\le j\le d.
\]
Let $S_0=0$ be  the original position and $S_n=\sum_{k=1}^nX_k$ be the position at time $n$ ($n\geq1$), where
\[
X_k\in\{\pm e_1,\ldots,\pm e_d\},\qquad
\Pp(X_1=e_j)=\Pp(X_1=-e_j)=\frac1{2d},\quad 1\le j\le d.
\]
For a memory parameter $\theta\in[0,1]$, set
\begin{equation}\label{eq:admissible}
a=\frac{2d\theta-1}{2d-1}\in
\mathcal A_d:=\left[-\frac1{2d-1},\frac12\right]\cap(-1,\infty).
\end{equation}
Let $\F_0$ be trivial and $\F_n=\sigma(X_1,\ldots,X_n)$. For $n\ge1$,
\begin{equation}\label{eq:transition}
\Pp(X_{n+1}=u\mid\F_n)
=\frac an\sum_{j=1}^n\ind_{\{X_j=u\}}+\frac{1-a}{2d},
\qquad u\in\{\pm e_1,\ldots,\pm e_d\}.
\end{equation}
The cases $a<1/2$ and $a=1/2$ correspond to the diffusive and critical regimes,
respectively. The endpoint $a=-1$ for $d=1$ is excluded throughout the paper.

For multidimensional martingales, Fang \textit{et al.}~\cite{FKS} established
normal approximation bounds in $p$-Wasserstein distance. Their Gaussian
convolution estimate with fixed covariance has linear dependence on $p$,
whereas the bound for general random predictable quadratic variation contains
an additional cube-root term. For the MERW, the conditional covariance depends
on both the coordinate occupation matrix and the second-order position tensor.
Their recursions allow a more precise treatment of this term.

In this paper, we establish a Gaussian convolution estimate with random
covariance under uniform boundedness and relative ellipticity conditions.
Combining this estimate with conditional expectations of the remaining
predictable covariance, we obtain $p$-Wasserstein bounds for the 
MERW. We further determine the covariance bias and derive radial Cram\'{e}r-type
moderate deviations by the results of Fang and Koike~\cite{FangKoike2023}.
The stated moderate deviation ranges are sufficient ranges; their optimality
is not asserted.

The paper is organized as follows. In Section~\ref{sec:main}, we present our
main results. In Section~\ref{sec:preliminaries}, we collect the Gaussian
convolution and covariance estimates. Section~\ref{sec:proofs} is devoted to
the proofs of the main results and a comparison with general martingale bounds.

Throughout the paper, $|x|$ denotes the Euclidean norm and $\norm A_{\op}$ the
operator norm. For a random vector $U$ and a random matrix $A$, write
\[
\norm U_q=(\E|U|^q)^{1/q},\qquad
\norm A_q=(\E\norm A_{\op}^q)^{1/q}.
\]
For probability distributions $\mu,\nu$ with finite $p$th moments, define
\[
W_{p,2}(\mu,\nu)
=\inf_{\Lw(U)=\mu,\,\Lw(V)=\nu}(\E|U-V|^p)^{1/p}.
\]
We reserve $p\ge1$ for the order of this distance. The notation $A\preceq B$
means that $B-A$ is positive semidefinite. We write $u_n\lesssim v_n$ if
$u_n\le Cv_n$, $u_n\asymp v_n$ if both inequalities hold, and $u_n\sim v_n$
if $u_n/v_n\to1$. Unless otherwise stated, $c,C$, possibly with indices,
are positive constants depending only on the fixed parameters $a,d$; their
values may vary from line to line. All asymptotic relations are as $n\to\infty$.
Auxiliary Gaussian vectors are defined on a suitable extension of the
probability space.

\section{Main results}\label{sec:main}
\quad
For $n\ge1$, set
\begin{equation}\label{eq:normalization}
b_n=\frac{\Gamma(a+1)\Gamma(n)}{\Gamma(n+a)},\qquad
v_n=\sum_{k=1}^nb_k^2,\qquad M_n=b_nS_n,\qquad
T_n=\frac{\sqrt d\,M_n}{\sqrt{v_n}}.
\end{equation}
By~\cite{BL} and Stirling's formula,
\begin{equation}\label{eq:weight_asymptotics}
b_n\sim\Gamma(a+1)n^{-a},\qquad
v_n\sim
\begin{cases}
\dfrac{\Gamma(a+1)^2}{1-2a}n^{1-2a},&a<1/2,\\[2mm]
\dfrac\pi4\log n,&a=1/2.
\end{cases}
\end{equation}
Define
\begin{equation}\label{eq:scales}
w_{n,k}=\frac{b_k^2}{v_n},\qquad
R_{n,k}=\sum_{\ell=k}^nw_{n,\ell},\qquad
L_n=\sum_{k=1}^n\frac{w_{n,k}^{3/2}}{R_{n,k}},
\end{equation}
with $R_{n,n+1}=0$, and write
\begin{equation}\label{eq:covariance_normalization}
\mathsf K_n=\Cov(T_n),\qquad
\widehat T_n=\mathsf K_n^{-1/2}T_n,\qquad
\delta_n=\frac{a^2}{v_n}\sum_{k=1}^{n-1}
\left(\frac{b_{k+1}}k\right)^2\E|S_k|^2,\qquad
\eta_n=L_n+\delta_n.
\end{equation}

We have the following Wasserstein bounds.
\begin{theorem}\label{thm:wasserstein}
Fix $d\ge1$ and $a\in\mathcal A_d$. Then
\begin{equation}\label{eq:exact_covariance}
\mathsf K_n=(1-\delta_n)I_d,\qquad
c_{a,d}I_d\preceq\mathsf K_n\preceq I_d.
\end{equation}
For all $n\ge1$ and $p\ge1$,
\begin{align}
\Wp\bigl(\Lw(T_n),N(0,\mathsf K_n)\bigr)
&\le C_{a,d}pL_n,\label{eq:main_covmatched}\\
\Wp\bigl(\Lw(\widehat T_n),N(0,I_d)\bigr)
&\le C_{a,d}pL_n,\label{eq:main_corrected}\\
\Wp\bigl(\Lw(T_n),N(0,I_d)\bigr)
&\le C_{a,d}(pL_n+\sqrt p\,\delta_n)
\le C_{a,d}p\eta_n.\label{eq:uncorrected_refined}
\end{align}
Moreover,
\begin{equation}\label{eq:Ln_asymp}
L_n\asymp
\begin{cases}
\dfrac{\log(en)}{\sqrt n},&a\le1/3,\\[2mm]
n^{-3(1/2-a)},&1/3<a<1/2,\\
(\log n)^{-3/2},&a=1/2,
\end{cases}
\end{equation}
and $\delta_n=0$ for $a=0$, whereas for $a\ne0$,
\begin{equation}\label{eq:delta_rates}
\delta_n\asymp
\begin{cases}
n^{-1},&a<0,\\
n^{-(1-2a)},&0<a<1/2,\\
(\log n)^{-1},&a=1/2.
\end{cases}
\end{equation}
Consequently,
\begin{equation}\label{eq:eta_rates}
\eta_n\asymp
\begin{cases}
\dfrac{\log(en)}{\sqrt n},&a\le1/4,\\[2mm]
n^{-(1-2a)},&1/4<a<1/2,\\
(\log n)^{-1},&a=1/2.
\end{cases}
\end{equation}
\end{theorem}

The covariance bias also gives the following lower bounds.
\begin{proposition}\label{prop:lower_bound}
For $n\ge1$ and $p\ge2$,
\begin{equation}\label{eq:lower_L2}
\Wp\bigl(\Lw(T_n),N(0,I_d)\bigr)
\ge\sqrt d\bigl(1-\sqrt{1-\delta_n}\bigr).
\end{equation}
If $a>1/4$, then, for all $p\ge1$ and all sufficiently large $n$,
\begin{equation}\label{eq:lower_all_p}
\Wp\bigl(\Lw(T_n),N(0,I_d)\bigr)\ge c_{a,d}\delta_n.
\end{equation}
\end{proposition}

Indeed, $L_n=o(\delta_n)$ for $a>1/4$, while at $a=1/4$,
$\delta_n\asymp n^{-1/2}=o(L_n)$.

We next give the radial Cram\'{e}r-type moderate deviations. Let
$Z\sim N(0,I_d)$.
\begin{theorem}\label{thm:cramer}
There exist $C,c>0$, depending only on $a,d$, such that, for all sufficiently
large $n$,
\begin{align}
\left|\frac{\Pp(|\widehat T_n|>x)}{\Pp(|Z|>x)}-1\right|
&\le C(1+x)(1+|\log L_n|+x^2)L_n,
&&0\le x\le cL_n^{-1/3},\label{eq:cramer_corrected}\\
\left|\frac{\Pp(|T_n|>x)}{\Pp(|Z|>x)}-1\right|
&\le C(1+x)(1+|\log\eta_n|+x^2)\eta_n,
&&0\le x\le c\eta_n^{-1/3}.\label{eq:cramer_uncorrected}
\end{align}
In particular,
\begin{align}
\sup_{0\le x\le x_n}
\left|\frac{\Pp(|\widehat T_n|>x)}{\Pp(|Z|>x)}-1\right|&\to0,
&&x_n=o(L_n^{-1/3}),\label{eq:cramer_corrected_uniform}\\
\sup_{0\le x\le x_n}
\left|\frac{\Pp(|T_n|>x)}{\Pp(|Z|>x)}-1\right|&\to0,
&&x_n=o(\eta_n^{-1/3}).\label{eq:cramer_uncorrected_uniform}
\end{align}
The respective ranges are determined by
\begin{equation}\label{eq:cramer_ranges}
\begin{aligned}
L_n^{-1/3}&\asymp
\begin{cases}
\dfrac{n^{1/6}}{(\log n)^{1/3}},&a\le1/3,\\[2mm]
n^{1/2-a},&1/3<a<1/2,\\
(\log n)^{1/2},&a=1/2,
\end{cases}\\[2mm]
\eta_n^{-1/3}&\asymp
\begin{cases}
\dfrac{n^{1/6}}{(\log n)^{1/3}},&a\le1/4,\\[2mm]
n^{(1-2a)/3},&1/4<a<1/2,\\
(\log n)^{1/3},&a=1/2.
\end{cases}
\end{aligned}
\end{equation}
\end{theorem}

\section{Preliminary lemmas}\label{sec:preliminaries}
\subsection{Gaussian convolution and replacement estimates}
\quad
We first extend the Gaussian convolution estimate of
Fang \textit{et al.}~\cite[Proposition 2.1]{FKS} to a random covariance.

\begin{lemma}\label{lem:local_random_cov}
Let $\mathcal G$ be a sub-$\sigma$-field, and suppose that
\[
\E(\xi\mid\mathcal G)=0,\quad
V=\E(\xi\xi^\top\mid\mathcal G),\quad
\bar A=\E(A\mid\mathcal G),\quad D=A-\bar A,\quad
\Sigma=V+\bar A\succ0,
\]
where $A$ is a random symmetric positive definite matrix. Conditionally on
$\mathcal G$, let $Z\sim N(0,I_d)$ be independent of $(\xi,A)$. Assume that
there exist fixed $c,C,K,\varepsilon_0>0$ and a $\mathcal G$-measurable
$0<\varepsilon\le\varepsilon_0$ such that the following relation holds almost surely,
\begin{equation}\label{eq:local_assumptions}
\begin{gathered}
c\Sigma\preceq A,\bar A\preceq C\Sigma,\qquad
|\Sigma^{-1/2}\xi|\le K\sqrt\varepsilon,\\
\norm{\Sigma^{-1/2}D\Sigma^{-1/2}}_{\op}\le K\varepsilon,
\qquad
\norm{\Sigma^{-1/2}V\Sigma^{-1/2}}_{\op}\le K\varepsilon.
\end{gathered}
\end{equation}
Then, for all $p\ge1$,
\begin{equation}\label{eq:local_main_bound}
\Wp\bigl(\Lw(\xi+A^{1/2}Z\mid\mathcal G),N(0,\Sigma)\bigr)
\le C_{c,C,K,\varepsilon_0,d}\norm{\Sigma^{1/2}}_{\op}
p\varepsilon^{3/2}.
\end{equation}
\end{lemma}

\begin{proof}
We work under a regular conditional law given $\mathcal G$. By applying
$\Sigma^{-1/2}$, it suffices to prove the assertion for $\Sigma=I_d$.
The normalized Gaussian term is represented by its covariance. Write
\[
\alpha=\norm{\xi}_\infty,\qquad
\delta=\norm D_\infty,\qquad v=\norm V_{\op};
\qquad \alpha\le K\sqrt\varepsilon,\quad
\delta,v\le K\varepsilon,\quad V+\bar A=I_d.
\]
In this proof, $C_*$ depends only on $c,C,K,\varepsilon_0,d$.

Set $W=\xi+A^{1/2}Z_1$, where $Z_1,Z_2$ are independent standard Gaussian
vectors independent of $(\xi,A)$. By~\cite[Lemma 5.1]{FKS},
\begin{equation}\label{eq:transport}
\Wp(\Lw(W),N(0,I_d))
\le2\int_0^\infty\norm{\rho_t(F_t)}_p\,dt,
\end{equation}
where, with $s=e^{-t}$,
\[
F_t=sW+\sqrt{1-s^2}Z_2,\qquad
\rho_t(F_t)=s\E\left[W-\frac{sZ_2}{\sqrt{1-s^2}}\,\middle|\,F_t\right].
\]
For $0\le u\le1$, define
\begin{equation}\label{eq:interpolation}
C_u=(1-s^2)I_d+s^2(\bar A+uD),\qquad
L_u=C_u^{1/2},\qquad G_u=s\xi+L_uZ,
\end{equation}
where $Z$ is a standard Gaussian vector independent of $(\xi,A)$.
Gaussian regression gives
\begin{equation}\label{eq:score_split}
\rho_t(G_1)=s m_1^A(G_1)+s^2m_1^D(G_1),
\end{equation}
with
\begin{equation}\label{eq:h_definitions}
\begin{aligned}
h_u^A(\xi,D,z)&=\xi-sVC_u^{-1/2}z,
&m_u^A(G_u)&=\E[h_u^A(\xi,D,Z)\mid G_u],\\
h_u^D(\xi,D,z)&=DC_u^{-1/2}z,
&m_u^D(G_u)&=\E[h_u^D(\xi,D,Z)\mid G_u].
\end{aligned}
\end{equation}
Here $F_t\overset{\mathrm d}=G_1$ and $I_d-A=V-D$.

By~\eqref{eq:local_assumptions},
\[
\min\{1,c\}I_d\preceq C_u\preceq\max\{1,C\}I_d.
\]
Differentiating $L_u^2=C_u$, we obtain
\[
L_u\dot L_u+\dot L_uL_u=s^2D,\qquad
\dot L_u=\int_0^\infty e^{-rL_u}(s^2D)e^{-rL_u}\,dr,
\]
and hence
\begin{equation}\label{eq:matrix_derivatives}
\norm{\dot L_u}_{\op}
+\norm{L_u^{-1}\dot L_u}_{\op}
+\norm{\partial_uC_u^{-1/2}}_{\op}
\le C_*s^2\delta.
\end{equation}

Let $\Theta=(\xi,D)$ and
\[
q_u(y\mid\Theta)
=\frac{\exp\{-\tfrac12(y-s\xi)^\top C_u^{-1}(y-s\xi)\}}
{(2\pi)^{d/2}(\det C_u)^{1/2}},\qquad
z_u=L_u^{-1}(y-s\xi).
\]
For either function $h_u$ in~\eqref{eq:h_definitions}, Bayes' formula gives
\[
m_u(y)=
\frac{\E[h_u(\Theta,z_u)q_u(y\mid\Theta)]}{\E[q_u(y\mid\Theta)]}.
\]
The boundedness of $\xi,D$ and uniform ellipticity of $C_u$ justify
differentiation under the expectations on compact sets; the denominator is
strictly positive. Therefore,
\begin{equation}\label{eq:conditional_derivatives}
\begin{aligned}
\partial_um_u(y)
&=\E[\partial_uh_u+(\partial_zh_u)\dot z_u\mid G_u=y]
+\Cov(h_u,\ell_u\mid G_u=y),\\
\nabla m_u(y)
&=\E[(\partial_zh_u)L_u^{-1}\mid G_u=y]
+\Cov(h_u,-L_u^{-1}z_u\mid G_u=y),
\end{aligned}
\end{equation}
where
\[
\ell_u=\partial_u\log q_u
=\frac12\left[z_u^\top L_u^{-1}(s^2D)L_u^{-1}z_u
-\tr\{L_u^{-1}(s^2D)L_u^{-1}\}\right],
\qquad \dot z_u=-L_u^{-1}\dot L_uz_u.
\]
At $y=G_u$, $z_u=Z$ and $\dot G_u=\dot L_uZ$. Thus, for $q\ge1$,
\begin{equation}\label{eq:derivative_norms}
\norm{\ell_u(G_u)}_q\le C_*q s^2\delta,\qquad
\norm{\dot z_u(G_u)}_q+\norm{\dot G_u}_q
\le C_*\sqrt q\,s^2\delta.
\end{equation}
Also, conditional Jensen's and H\"older's inequalities yield
\begin{equation}\label{eq:conditional_cov_bound}
\norm{\Cov(H,R\mid G)}_p\le2\norm H_{2p}\norm R_{2p}.
\end{equation}

By~\eqref{eq:matrix_derivatives} and~\eqref{eq:h_definitions},
\[
\begin{aligned}
\norm{h_u^A}_q&\le\alpha+C_*sv\sqrt q,
&\norm{\partial_zh_u^A}_{\op}&\le C_*sv,
&\norm{\partial_uh_u^A}_q&\le C_*\sqrt q\,s^3v\delta,\\
\norm{h_u^D}_q&\le C_*\sqrt q\,\delta,
&\norm{\partial_zh_u^D}_{\op}&\le C_*\delta,
&\norm{\partial_uh_u^D}_q&\le C_*\sqrt q\,s^2\delta^2.
\end{aligned}
\]
Combining~\eqref{eq:conditional_derivatives}--\eqref{eq:conditional_cov_bound},
we obtain
\[
\norm{\nabla m_u^A(G_u)}_{2p}\le C_*(\sqrt p\,\alpha+psv),\qquad
\norm{\nabla m_u^D(G_u)}_{2p}\le C_*p\delta,
\]
and
\begin{align*}
\norm{\frac d{du}m_u^A(G_u)}_p
&\le\norm{\partial_um_u^A(G_u)}_p
+\norm{\nabla m_u^A(G_u)\dot G_u}_p
\le C_*(ps^2\alpha\delta+p^{3/2}s^3v\delta),\\
\norm{\frac d{du}m_u^D(G_u)}_p
&\le\norm{\partial_um_u^D(G_u)}_p
+\norm{\nabla m_u^D(G_u)\dot G_u}_p
\le C_*p^{3/2}s^2\delta^2.
\end{align*}
These bounds give absolute continuity in $L^p$. Integrating over $u\in[0,1]$
yields
\begin{equation}\label{eq:interpolation_bounds}
\begin{aligned}
\norm{m_1^A(G_1)-m_0^A(G_0)}_p
&\le C_*(ps^2\alpha\delta+p^{3/2}s^3v\delta),\\
\norm{m_1^D(G_1)-m_0^D(G_0)}_p
&\le C_*p^{3/2}s^2\delta^2.
\end{aligned}
\end{equation}

Since $C_0=I_d-s^2V$ is deterministic and $\E D=0$,
\[
m_0^D(G_0)=\E[DC_0^{-1/2}Z\mid Z+sC_0^{-1/2}\xi],\qquad
\E[DC_0^{-1/2}Z\mid Z]=0.
\]
Apply~\cite[Lemma 5.3]{FKS} with
\[
X=0,\qquad Y=sC_0^{-1/2}\xi,\qquad H=DC_0^{-1/2}Z,
\]
and $D$ as the auxiliary random element. Since
\[
\norm H_q\le C_*\sqrt q\,\delta,\qquad
\norm{\partial_ZH}_{L^\infty(\op)}\le C_*\delta,\qquad
\norm Y_\infty\le C_*s\alpha,
\]
the choices $(a,c)=(\infty,p)$, $(\alpha_0,\beta_0)=(2p,2p)$ and
$(\alpha_2,\beta_2,\gamma_2)=(p,\infty,\infty)$ in that lemma give
\begin{equation}\label{eq:centering_bound}
\norm{m_0^D(G_0)}_p
\le C_*(ps\delta\alpha+\sqrt p\,s^2\delta\alpha^2).
\end{equation}
For the fixed-covariance term $\rho_t^{(0)}=s m_0^A(G_0)$, the score estimate
in the proof of~\cite[Proposition 2.1]{FKS} gives
\[
2\int_0^\infty\norm{\rho_t^{(0)}}_p\,dt
\le C_*p\norm\xi_{3p}^3\le C_*p\alpha^3.
\]
Moreover, by~\eqref{eq:score_split},
\[
\rho_t-\rho_t^{(0)}
=s[m_1^A(G_1)-m_0^A(G_0)]
+s^2[m_1^D(G_1)-m_0^D(G_0)]+s^2m_0^D(G_0).
\]
Using~\eqref{eq:transport}, \eqref{eq:interpolation_bounds}
and~\eqref{eq:centering_bound}, and integrating $s^j=e^{-jt}$, we obtain
\begin{equation}\label{eq:refined_local_bound}
\Wp(\Lw(W),N(0,I_d))
\le C_*\{p\alpha^3+p\alpha\delta+\sqrt p\,\delta\alpha^2
+p^{3/2}(v\delta+\delta^2)\}.
\end{equation}

On the other hand, the Sylvester identity yields
\[
A^{1/2}-\bar A^{1/2}
=\int_0^\infty e^{-rA^{1/2}}D e^{-r\bar A^{1/2}}\,dr,
\qquad
\norm{A^{1/2}-\bar A^{1/2}}_{\op}\le\frac{\delta}{2\sqrt c}.
\]
Coupling with the same $(\xi,Z)$ and applying~\cite[Proposition 2.1]{FKS},
\[
\Wp(\Lw(W),N(0,I_d))\le C_*(p\alpha^3+\sqrt p\,\delta).
\]
Together with~\eqref{eq:refined_local_bound}, this gives
\begin{align*}
\Wp(\Lw(W),N(0,I_d))
&\le C_*\{p\varepsilon^{3/2}
+\min(\sqrt p\,\varepsilon,p^{3/2}\varepsilon^2)\}\\
&\le C_*p\varepsilon^{3/2},
\end{align*}
where $\min(x,y)\le\sqrt{xy}$ was used in the last inequality.
Applying $\Sigma^{1/2}$ proves~\eqref{eq:local_main_bound}.
\end{proof}

We next use conditional expectations of the remaining predictable covariance.
\begin{lemma}\label{lem:replacement}
Let $(\xi_k,\F_k)_{k=1}^n$ be $\R^d$-valued martingale differences with
$\E|\xi_k|^{\max\{2,p\}}<\infty$, where $p\ge1$ and $\F_0$ is trivial.
For a deterministic matrix $M\succ0$, set
\begin{equation}\label{eq:remaining_covariance}
\begin{gathered}
U_k=\sum_{i=1}^k\xi_i,\qquad
V_k=\E(\xi_k\xi_k^\top\mid\F_{k-1}),\qquad
\bar\Sigma=\E\sum_{i=1}^nV_i,\\
Q_k=M+\E\left[\sum_{i=k+1}^nV_i\,\middle|\,\F_k\right],\qquad
\Gamma_k=Q_k-\E(Q_k\mid\F_{k-1}),\qquad 0\le k\le n,
\end{gathered}
\end{equation}
where $U_0=0$ and $\Gamma_k$ is defined for $k\ge1$. Then
\begin{equation}\label{eq:Q_recursion}
Q_{k-1}=V_k+\E(Q_k\mid\F_{k-1}),\qquad
Q_0=M+\bar\Sigma,\qquad Q_n=M.
\end{equation}
For $Z\sim N(0,I_d)$ independent of $\F_n$,
\begin{equation}\label{eq:global_telescoping}
\begin{aligned}
&\Wp\bigl(\Lw(U_n+M^{1/2}Z),N(0,\bar\Sigma+M)\bigr)\\
&\quad\le\sum_{k=1}^n
\left\{\E\left[
\Wp^p\bigl(\Lw(\xi_k+Q_k^{1/2}Z\mid\F_{k-1}),N(0,Q_{k-1})\bigr)
\right]\right\}^{1/p}.
\end{aligned}
\end{equation}
\end{lemma}

\begin{proof}
The tower property gives
\[
V_k+\E(Q_k\mid\F_{k-1})
=M+\E\left[\sum_{i=k}^nV_i\,\middle|\,\F_{k-1}\right]=Q_{k-1}.
\]
In particular,
\begin{equation}\label{eq:local_cov_matching}
\E[(\xi_k+Q_k^{1/2}Z)(\xi_k+Q_k^{1/2}Z)^\top\mid\F_{k-1}]
=V_k+\E(Q_k\mid\F_{k-1})=Q_{k-1}.
\end{equation}
Let $Z_0,\ldots,Z_n$ be i.i.d. standard Gaussian vectors independent of $\F_n$,
and set $H_k=U_k+Q_k^{1/2}Z_k$. Conditional Jensen's inequality gives
$\E\norm{Q_k}_{\op}^{p/2}<\infty$ for $p\ge2$; square integrability suffices
for $1\le p<2$. Thus all the following distances are finite. Since
\[
H_0\sim N(0,\bar\Sigma+M),\qquad
H_n\overset{\mathrm d}=U_n+M^{1/2}Z,
\]
the triangle inequality and conditional coupling imply
\begin{align*}
\Wp(\Lw(H_n),\Lw(H_0))
&\le\sum_{k=1}^n\Wp(\Lw(H_k),\Lw(H_{k-1}))\\
&\le\sum_{k=1}^n
\left\{\E\left[
\Wp^p\bigl(\Lw(\xi_k+Q_k^{1/2}Z\mid\F_{k-1}),N(0,Q_{k-1})\bigr)
\right]\right\}^{1/p}.
\end{align*}
This proves~\eqref{eq:global_telescoping}; see also the conditional coupling
argument in~\cite[Theorem 2.1]{FKS}.
\end{proof}

The following consequence of Fang and Koike~\cite{FangKoike2023} will be used
for moderate deviations.
\begin{lemma}\label{lem:cramer_transfer}
Fix $d\ge1$ and $A>0$. Let $Z\sim N(0,I_d)$ and suppose that
\[
\Wp(\Lw(W),\Lw(Z))\le Ap\Delta,\qquad p\ge1.
\]
There exist $C,c,\Delta_0>0$, depending only on $A,d$, such that, for
$0<\Delta\le\Delta_0$ and $0\le x\le c\Delta^{-1/3}$,
\begin{equation}\label{eq:generic_cramer_transfer}
\left|\frac{\Pp(|W|>x)}{\Pp(|Z|>x)}-1\right|
\le C(1+x)(1+|\log\Delta|+x^2)\Delta.
\end{equation}
\end{lemma}

\begin{proof}
For $d\ge2$, apply~\cite[Theorem 4.2]{FangKoike2023} with growth exponent
$\alpha=1$ and $p_0=\Delta^{-2/3}$. With
$\kappa(d)=2^{d/2-1}\Gamma(d/2)$, for sufficiently small $\Delta$,
\[
|\log\Delta|\vee\log\kappa(d)\le p_0/4,\qquad
d(d\log d)\Delta\le B_1,\qquad
\min\{\Delta^{-1/3},\sqrt{p_0}\}=\Delta^{-1/3},
\]
where $B_1>0$ is fixed. The theorem gives
\[
\left|\frac{\Pp(|W|>x)}{\Pp(|Z|>x)}-1\right|
\le C(1+x)(|\log\Delta|+d\log d+x^2)\Delta,
\]
which implies~\eqref{eq:generic_cramer_transfer}. For $d=1$, apply
\cite[Theorem 2.2]{FangKoike2023} to $W$ and $-W$ and use
$\Pp(|Z|>x)=2\Pp(Z>x)$.
\end{proof}

\subsection{Covariance estimates for the MERW}

Set
\[
\mathcal O_m=\sum_{j=1}^mX_jX_j^\top,\qquad
D_m^{\rm occ}=\mathcal O_m-\frac mdI_d,\qquad
\varepsilon_1=X_1,\qquad
\varepsilon_{m+1}=X_{m+1}-\frac amS_m.
\]
By~\eqref{eq:transition}, or~\cite[equations (2.3) and (A.5)]{BL},
\begin{equation}\label{eq:model_moments}
\E(X_{m+1}\mid\F_m)=\frac amS_m,\qquad
\E(X_{m+1}X_{m+1}^\top\mid\F_m)
=\frac am\mathcal O_m+\frac{1-a}{d}I_d.
\end{equation}
Consequently, with $C_1=I_d/d$,
\begin{equation}\label{eq:innovationcov}
C_{m+1}:=\E(\varepsilon_{m+1}\varepsilon_{m+1}^\top\mid\F_m)
=\frac{1-a}{d}I_d+\frac am\mathcal O_m-\frac{a^2}{m^2}S_mS_m^\top.
\end{equation}
Since $b_{m+1}(1+a/m)=b_m$,
\begin{equation}\label{eq:martingale_array}
M_n=\sum_{k=1}^nb_k\varepsilon_k,\qquad
T_n=\sum_{k=1}^n\xi_{n,k},\qquad
\xi_{n,k}=\sqrt{d w_{n,k}}\,\varepsilon_k,\qquad
|\xi_{n,k}|\le2\sqrt{d w_{n,k}}.
\end{equation}
For $k\ge2$,
\begin{equation}\label{eq:Vdecomp}
\begin{aligned}
V_{n,k}&:=\E(\xi_{n,k}\xi_{n,k}^\top\mid\F_{k-1})
=w_{n,k}I_d+A_{n,k}^{\rm occ}+B_{n,k}^{S},\\
A_{n,k}^{\rm occ}&=\frac{da b_k^2}{v_n(k-1)}D_{k-1}^{\rm occ},
\qquad
B_{n,k}^{S}=-\frac{da^2b_k^2}{v_n(k-1)^2}S_{k-1}S_{k-1}^\top.
\end{aligned}
\end{equation}
For $k=1$, set $A_{n,1}^{\rm occ}=B_{n,1}^{S}=0$, so that
$V_{n,1}=w_{n,1}I_d$.

For fixed $n$, use~\eqref{eq:remaining_covariance} with $(\xi_{n,k})_{k=1}^n$,
and write
\begin{equation}\label{eq:Gamma_decomp}
\begin{aligned}
\Gamma_{n,k}
&=Q_k-\E(Q_k\mid\F_{k-1})\\
&=\sum_{i=k+1}^n
\{\E(V_{n,i}\mid\F_k)-\E(V_{n,i}\mid\F_{k-1})\}
=\Gamma_{n,k}^{\rm occ}+\Gamma_{n,k}^{S},
\end{aligned}
\end{equation}
where the two terms correspond to $A_{n,i}^{\rm occ}$ and $B_{n,i}^{S}$
in~\eqref{eq:Vdecomp}.

\begin{lemma}\label{lem:ellipticity}
For fixed $a\in\mathcal A_d$,
\begin{equation}\label{eq:V_uniform_ellipticity}
c_{a,d}I_d\preceq C_k\preceq I_d,\qquad
c_{a,d}w_{n,k}I_d\preceq V_{n,k}\preceq dw_{n,k}I_d.
\end{equation}
If $M=w_{n,n}I_d$ and $r_{n,k}=R_{n,k}+w_{n,n}$, then
\begin{equation}\label{eq:Q_scale}
c_{a,d}r_{n,k}I_d\preceq
Q_{k-1},\ Q_k,\ \E(Q_k\mid\F_{k-1})
\preceq C_{a,d}r_{n,k}I_d.
\end{equation}
\end{lemma}

\begin{proof}
By the empirical covariance inequality,
\[
\frac{S_mS_m^\top}{m^2}\preceq\frac{\mathcal O_m}{m}\preceq I_d.
\]
Thus~\eqref{eq:innovationcov} gives
\[
C_{m+1}\succeq\frac{1-a}{d}I_d
+a(1-a)\frac{\mathcal O_m}{m}
\succeq
\begin{cases}
\dfrac{1-a}{d}I_d,&a\ge0,\\[1mm]
(1-a)(d^{-1}+a)I_d,&a<0.
\end{cases}
\]
Both constants are positive by~\eqref{eq:admissible}. Also,
\[
\tr C_{m+1}=1-\frac{a^2}{m^2}|S_m|^2\le1,
\]
so $C_{m+1}\preceq I_d$. Together with $C_1=I_d/d$ and
$V_{n,k}=d w_{n,k}C_k$, this proves~\eqref{eq:V_uniform_ellipticity}.

It follows that
\[
c_{a,d}(R_{n,k+1}+w_{n,n})I_d\preceq Q_k
\preceq C_{a,d}(R_{n,k+1}+w_{n,n})I_d.
\]
For $k<n$,
\[
\frac{w_{n,k}}{w_{n,k+1}}=\left(1+\frac ak\right)^2\le C_a,
\qquad R_{n,k+1}+w_{n,n}\asymp R_{n,k}+w_{n,n},
\]
and $r_{n,n}=2w_{n,n}$. Taking conditional expectations and using
\eqref{eq:Q_recursion} proves~\eqref{eq:Q_scale}.
\end{proof}

\begin{lemma}\label{lem:covariance_increments}
For $a<0$, set
\[
q_{n,k}=\frac{k^{-a}(n^{-a}-k^{-a})}{n^{1-2a}},
\]
and define
\begin{equation}\label{eq:hnk_def}
h_{n,k}=
\begin{cases}
w_{n,k},&a\ge0,\\
w_{n,k}+q_{n,k},&a<0.
\end{cases}
\end{equation}
Then
\begin{equation}\label{eq:Gamma_total_bound}
\norm{\Gamma_{n,k}^{S}}_{\op}\le C_{a,d}w_{n,k},\qquad
\norm{\Gamma_{n,k}}_{\op}\le C_{a,d}h_{n,k}.
\end{equation}
Moreover, for $a<0$,
\begin{equation}\label{eq:q_bounds}
q_{n,k}\le\frac Cn,\qquad
q_{n,k}\le C\frac{n-k+1}{n^2}\quad(k>n/2).
\end{equation}
\end{lemma}

\begin{proof}
Set $\Delta_k=X_kX_k^\top-\E(X_kX_k^\top\mid\F_{k-1})$. Then
\[
\E(\Delta_k\mid\F_{k-1})=0,\qquad
\norm{\Delta_k}_{\op}\le1,\qquad \Delta_1=D_1^{\rm occ}.
\]
By~\eqref{eq:model_moments},
\[
D_{m+1}^{\rm occ}=\left(1+\frac am\right)D_m^{\rm occ}+\Delta_{m+1},
\qquad b_mD_m^{\rm occ}=\sum_{j=1}^mb_j\Delta_j.
\]
Hence, for $r\ge k$,
\begin{equation}\label{eq:occ_Doob_increment}
\E(D_r^{\rm occ}\mid\F_k)-\E(D_r^{\rm occ}\mid\F_{k-1})
=\frac{b_k}{b_r}\Delta_k.
\end{equation}
Using $b_i^2/((i-1)b_{i-1})=b_i/(i-1+a)$ in~\eqref{eq:Vdecomp},
\begin{equation}\label{eq:Gamma_occ_exact}
\Gamma_{n,k}^{\rm occ}
=\frac{da b_k}{v_n}
\left(\sum_{i=k+1}^n\frac{b_i}{i-1+a}\right)\Delta_k.
\end{equation}
For $a>0$,
\[
\norm{\Gamma_{n,k}^{\rm occ}}_{\op}
\le\frac{Cb_k}{v_n}\sum_{i>k}i^{-a-1}
\le\frac{Cb_k k^{-a}}{v_n}\le Cw_{n,k},
\]
whereas $\Gamma_{n,k}^{\rm occ}=0$ for $a=0$. For $a<0$,
\[
\norm{\Gamma_{n,k}^{\rm occ}}_{\op}
\le\frac{Ck^{-a}}{n^{1-2a}}\sum_{i=k+1}^ni^{-a-1}
\le Cq_{n,k}.
\]
In this case,
\[
q_{n,k}=\frac1n\left(\frac kn\right)^{-a}
\left[1-\left(\frac kn\right)^{-a}\right]\le\frac1{4n},
\]
and, for $k>n/2$, the mean value theorem gives
$n^{-a}-k^{-a}\le C_an^{-a-1}(n-k)$, proving~\eqref{eq:q_bounds}.

For the second-order term, set
\[
\mathcal Q_m=S_mS_m^\top-\E(S_mS_m^\top),\qquad
\Xi_k=S_kS_k^\top-\E(S_kS_k^\top\mid\F_{k-1}).
\]
Then $\Xi_1=\mathcal Q_1$ and
\begin{equation}\label{eq:tensor_recursion}
\begin{aligned}
\mathcal Q_{m+1}
&=\left(1+\frac{2a}{m}\right)\mathcal Q_m
+\frac amD_m^{\rm occ}+\Xi_{m+1},\\
\Xi_{m+1}
&=S_m\varepsilon_{m+1}^\top+\varepsilon_{m+1}S_m^\top+\Delta_{m+1},
\qquad \norm{\Xi_k}_{\op}\le Ck.
\end{aligned}
\end{equation}
Define
\[
G_{m,k}=\prod_{\ell=k}^{m-1}\left(1+\frac{2a}{\ell}\right),
\qquad G_{k,k}=1.
\]
Since $\log(1+2a/\ell)=2a/\ell+O_a(\ell^{-2})$ for all sufficiently large
$\ell$, with the finitely many initial factors absorbed into $C_a$,
\begin{equation}\label{eq:Gmk_bound}
|G_{m,k}|\le C_a(m/k)^{2a},\qquad m\ge k\ge1.
\end{equation}
The estimate also holds when an initial factor vanishes.
Iterating~\eqref{eq:tensor_recursion} and using~\eqref{eq:occ_Doob_increment},
\begin{equation}\label{eq:tensor_increment}
\begin{aligned}
&\E(\mathcal Q_m\mid\F_k)-\E(\mathcal Q_m\mid\F_{k-1})\\
&\quad=G_{m,k}\Xi_k
+a b_k\Delta_k\sum_{\ell=k}^{m-1}\frac{G_{m,\ell+1}}{\ell b_\ell}.
\end{aligned}
\end{equation}
Therefore,
\[
\Gamma_{n,k}^{S}
=-\frac{da^2}{v_n}\sum_{i=k+1}^n\frac{b_i^2}{(i-1)^2}
\left(G_{i-1,k}\Xi_k
+a b_k\Delta_k\sum_{\ell=k}^{i-2}\frac{G_{i-1,\ell+1}}{\ell b_\ell}\right).
\]
By~\eqref{eq:Gmk_bound} and $b_i\asymp i^{-a}$,
\begin{align*}
\frac{\norm{\Xi_k}_{\op}}{v_n}
\sum_{i>k}\frac{b_i^2}{(i-1)^2}|G_{i-1,k}|
&\le\frac{Ck^{1-2a}}{v_n}\sum_{i>k}i^{-2}
\le Cw_{n,k},\\
\frac{b_k}{v_n}\sum_{i=k+1}^n\frac{b_i^2}{(i-1)^2}
\sum_{\ell=k}^{i-2}\frac{|G_{i-1,\ell+1}|}{\ell b_\ell}
&\le\frac{Cb_k}{v_n}\sum_{\ell\ge k}\frac{\ell^{-2a}}{\ell b_\ell}
\sum_{i>\ell}i^{-2}\\
&\le\frac{Cb_k}{v_n}\sum_{\ell\ge k}\ell^{-a-2}
\le\frac{Ck^{-1-2a}}{v_n}\le Cw_{n,k}.
\end{align*}
Here $a>-1$ ensures convergence of the last series. Combining these estimates
with~\eqref{eq:Gamma_occ_exact} proves~\eqref{eq:Gamma_total_bound}.
\end{proof}

\begin{lemma}\label{lem:scales_bias}
The covariance identity~\eqref{eq:exact_covariance} and the estimates
\eqref{eq:Ln_asymp}--\eqref{eq:eta_rates} hold.
\end{lemma}

\begin{proof}
For $a<1/2$, \eqref{eq:weight_asymptotics} gives
\[
w_{n,k}\asymp\frac{k^{-2a}}{n^{1-2a}},\qquad
R_{n,k}\asymp
\begin{cases}
1,&k\le n/2,\\
(n-k+1)/n,&k>n/2.
\end{cases}
\]
Thus
\begin{align*}
\sum_{k\le n/2}\frac{w_{n,k}^{3/2}}{R_{n,k}}
&\asymp n^{-3(1/2-a)}\sum_{k\le n/2}k^{-3a}
\asymp
\begin{cases}
n^{-1/2},&a<1/3,\\
n^{-1/2}\log n,&a=1/3,\\
n^{-3(1/2-a)},&a>1/3,
\end{cases}\\
\sum_{k>n/2}\frac{w_{n,k}^{3/2}}{R_{n,k}}
&\asymp\frac1{\sqrt n}\sum_{m\le\lceil n/2\rceil}\frac1m
\asymp\frac{\log n}{\sqrt n}.
\end{align*}
For $a=1/2$,
\[
w_{n,k}\asymp\frac1{k\log n},\qquad
R_{n,k}\asymp\frac1{\log n}\sum_{\ell=k}^n\frac1\ell,
\]
and hence
\begin{align*}
\sum_{k\le\sqrt n}\frac{w_{n,k}^{3/2}}{R_{n,k}}
&\asymp(\log n)^{-3/2}\sum_{k\le\sqrt n}k^{-3/2}
\asymp(\log n)^{-3/2},\\
\sum_{\sqrt n<k\le n/2}\frac{w_{n,k}^{3/2}}{R_{n,k}}
&\lesssim(\log n)^{-1/2}\sum_{k>\sqrt n}k^{-3/2}
\lesssim n^{-1/4}(\log n)^{-1/2},\\
\sum_{k>n/2}\frac{w_{n,k}^{3/2}}{R_{n,k}}
&\lesssim\frac1{\sqrt{n\log n}}
\sum_{m\le\lceil n/2\rceil}\frac1m
\lesssim\sqrt{\frac{\log n}{n}}.
\end{align*}
The last two bounds are $o((\log n)^{-3/2})$, proving~\eqref{eq:Ln_asymp}.

By coordinate symmetry and~\eqref{eq:Vdecomp},
\[
\E\mathcal O_k=\frac kdI_d,\qquad
\E(S_kS_k^\top)=\frac{\E|S_k|^2}{d}I_d,
\]
so that
\begin{align*}
\mathsf K_n=\E\sum_{k=1}^nV_{n,k}
&=I_d-\frac{da^2}{v_n}\sum_{k=1}^{n-1}
\left(\frac{b_{k+1}}k\right)^2\E(S_kS_k^\top)\\
&=(1-\delta_n)I_d.
\end{align*}
The lower bound in~\eqref{eq:exact_covariance} follows from
\eqref{eq:V_uniform_ellipticity} and $\sum_kw_{n,k}=1$; the upper bound follows
from $\delta_n\ge0$.

For $u_k=\E|S_k|^2$, \eqref{eq:model_moments} gives
\[
u_1=1,\qquad u_{k+1}=\left(1+\frac{2a}{k}\right)u_k+1.
\]
For $a<1/2$,
\[
\frac{u_{k+1}}{k+1}-\frac1{1-2a}
=\left(1-\frac{1-2a}{k+1}\right)
\left(\frac{u_k}{k}-\frac1{1-2a}\right)
\longrightarrow0,
\]
since the product of the factors tends to zero. For $a=1/2$,
\[
\frac{u_{k+1}}{k+1}=\frac{u_k}{k}+\frac1{k+1},\qquad
u_k=k\sum_{j=1}^k\frac1j.
\]
Consequently,
\[
u_k\sim\frac{k}{1-2a}\quad(a<1/2),\qquad
u_k\asymp k\log(ek)\quad(a=1/2).
\]
For $a\ne0$, substituting into~\eqref{eq:covariance_normalization} yields
\begin{equation*}
\delta_n\asymp
\begin{cases}
\displaystyle\frac1{n^{1-2a}}\sum_{k<n}k^{-1-2a}
\asymp n^{-1},&a<0,\\[2mm]
\displaystyle\frac1{n^{1-2a}}\sum_{k<n}k^{-1-2a}
\asymp n^{-(1-2a)},&0<a<1/2,\\[2mm]
\displaystyle\frac1{\log n}\sum_{k<n}\frac{\log(ek)}{k^2}
\asymp(\log n)^{-1},&a=1/2.
\end{cases}
\end{equation*}
For $a=0$, $\delta_n=0$. This proves~\eqref{eq:delta_rates}, and
\eqref{eq:eta_rates} follows from $\eta_n=L_n+\delta_n$.
\end{proof}

\section{Proofs of the main results}\label{sec:proofs}

\subsection{Proof of Theorem~\ref{thm:wasserstein}}
Apply Lemma~\ref{lem:replacement} to $(\xi_{n,k})_{k=1}^n$ with
$M=w_{n,n}I_d$. For the $k$th replacement, take
\[
\Sigma=Q_{k-1},\quad A=Q_k,\quad
\bar A=\E(Q_k\mid\F_{k-1}),\quad
D=\Gamma_{n,k},\quad V=V_{n,k},\quad
\varepsilon=\frac{h_{n,k}}{r_{n,k}}.
\]
By~\eqref{eq:Q_recursion}, $V+\bar A=\Sigma$. Lemmas~\ref{lem:ellipticity}
and~\ref{lem:covariance_increments} give
\begin{gather*}
c_{a,d}\Sigma\preceq A,\bar A\preceq C_{a,d}\Sigma,
\qquad \norm{\Sigma^{1/2}}_{\op}\le C\sqrt{r_{n,k}},\\
|\Sigma^{-1/2}\xi_{n,k}|\le C\sqrt{\frac{w_{n,k}}{r_{n,k}}},\qquad
\norm{\Sigma^{-1/2}V_{n,k}\Sigma^{-1/2}}_{\op}
\le C\frac{w_{n,k}}{r_{n,k}},\\
\norm{\Sigma^{-1/2}\Gamma_{n,k}\Sigma^{-1/2}}_{\op}
\le C\frac{h_{n,k}}{r_{n,k}}.
\end{gather*}
Since $w_{n,k}\le r_{n,k}$ and, for $a<0$,
$q_{n,k}\lesssim n^{-1}\asymp w_{n,n}\le r_{n,k}$,
\[
0<\frac{h_{n,k}}{r_{n,k}}\le C_{a,d}.
\]
Lemma~\ref{lem:local_random_cov} therefore yields
\begin{equation}\label{eq:local_swap_error}
\Wp\bigl(\Lw(\xi_{n,k}+Q_k^{1/2}Z\mid\F_{k-1}),N(0,Q_{k-1})\bigr)
\le C_{a,d}p\frac{h_{n,k}^{3/2}}{r_{n,k}}.
\end{equation}

For $a\ge0$,
\[
\sum_{k=1}^n\frac{h_{n,k}^{3/2}}{r_{n,k}}
\le\sum_{k=1}^n\frac{w_{n,k}^{3/2}}{R_{n,k}}=L_n.
\]
For $a<0$, by~\eqref{eq:q_bounds},
\begin{align*}
\sum_{k\le n/2}\frac{q_{n,k}^{3/2}}{r_{n,k}}
&\le Cn\cdot n^{-3/2}\le Cn^{-1/2},\\
\sum_{k>n/2}\frac{q_{n,k}^{3/2}}{r_{n,k}}
&\le C\sum_{m\le\lceil n/2\rceil}
\frac{(m/n^2)^{3/2}}{m/n}
\le\frac C{n^2}\sum_{m\le\lceil n/2\rceil}\sqrt m
\le Cn^{-1/2}.
\end{align*}
Together with $(x+y)^{3/2}\le C(x^{3/2}+y^{3/2})$ and~\eqref{eq:Ln_asymp},
\[
\sum_{k=1}^n\frac{h_{n,k}^{3/2}}{r_{n,k}}\le C_{a,d}L_n.
\]
It follows from~\eqref{eq:global_telescoping} and~\eqref{eq:local_swap_error}
that
\begin{equation}\label{eq:smoothed_global_bound}
\Wp\bigl(\Lw(T_n+M^{1/2}Z),N(0,\mathsf K_n+M)\bigr)
\le C_{a,d}pL_n.
\end{equation}
Moreover,
\[
\norm{M^{1/2}Z}_p\le C_d\sqrt p\sqrt{w_{n,n}},\qquad
L_n\ge\frac{w_{n,n}^{3/2}}{R_{n,n}}=\sqrt{w_{n,n}}.
\]
Removing the Gaussian smoothing from both distributions in
\eqref{eq:smoothed_global_bound} gives
\[
\Wp\bigl(\Lw(T_n),N(0,\mathsf K_n)\bigr)
\le C_{a,d}pL_n+2C_d\sqrt p\sqrt{w_{n,n}}
\le C_{a,d}pL_n.
\]
This proves~\eqref{eq:main_covmatched}. By~\eqref{eq:exact_covariance},
\[
\Wp\bigl(\Lw(\widehat T_n),N(0,I_d)\bigr)
\le\norm{\mathsf K_n^{-1/2}}_{\op}
\Wp\bigl(\Lw(T_n),N(0,\mathsf K_n)\bigr)
\le C_{a,d}pL_n.
\]
Finally, the Gaussian scaling coupling gives
\begin{equation}\label{eq:gaussian_cov_compare}
\begin{aligned}
\Wp\bigl(N(0,\mathsf K_n),N(0,I_d)\bigr)
&\le\norm Z_p\bigl(1-\sqrt{1-\delta_n}\bigr)\\
&\le C_d\sqrt p\,\delta_n.
\end{aligned}
\end{equation}
Combining this with~\eqref{eq:main_covmatched} proves
\eqref{eq:uncorrected_refined}. The remaining assertions follow from
Lemma~\ref{lem:scales_bias}.

\subsection{Proof of Proposition~\ref{prop:lower_bound}}
For every coupling $(U,Z)$ with $\Lw(U)=\Lw(T_n)$ and $Z\sim N(0,I_d)$,
\[
\norm{U-Z}_2\ge\norm Z_2-\norm U_2
=\sqrt d\bigl(1-\sqrt{1-\delta_n}\bigr).
\]
Taking the infimum and using $W_{p,2}\ge W_{2,2}$ for $p\ge2$
proves~\eqref{eq:lower_L2}. Also, the scaling coupling and the
$1$-Lipschitz function $x\mapsto|x|$ give
\[
W_{1,2}\bigl(N(0,\mathsf K_n),N(0,I_d)\bigr)
=\E|Z|\bigl(1-\sqrt{1-\delta_n}\bigr).
\]
For $a>1/4$, $L_n=o(\delta_n)$, and hence, by~\eqref{eq:main_covmatched},
\begin{align*}
\Wp\bigl(\Lw(T_n),N(0,I_d)\bigr)
&\ge W_{1,2}\bigl(\Lw(T_n),N(0,I_d)\bigr)\\
&\ge\E|Z|\bigl(1-\sqrt{1-\delta_n}\bigr)-C_{a,d}L_n\\
&\ge\tfrac12\E|Z|\delta_n-C_{a,d}L_n
\ge c_{a,d}\delta_n
\end{align*}
for all sufficiently large $n$. This proves~\eqref{eq:lower_all_p}.

\subsection{Proof of Theorem~\ref{thm:cramer}}
By~\eqref{eq:main_corrected} and~\eqref{eq:uncorrected_refined},
\[
\Wp\bigl(\Lw(\widehat T_n),N(0,I_d)\bigr)\le C_{a,d}pL_n,
\qquad
\Wp\bigl(\Lw(T_n),N(0,I_d)\bigr)\le C_{a,d}p\eta_n.
\]
Applying Lemma~\ref{lem:cramer_transfer} with $\Delta=L_n$ and
$\Delta=\eta_n$, respectively, proves~\eqref{eq:cramer_corrected}
and~\eqref{eq:cramer_uncorrected}. Since
\[
x_n=o(\Delta^{-1/3})
\quad\Longrightarrow\quad
\Delta(1+x_n)(1+|\log\Delta|+x_n^2)=o(1),
\]
the uniform limits follow. The explicit ranges in~\eqref{eq:cramer_ranges}
are immediate from~\eqref{eq:Ln_asymp} and~\eqref{eq:eta_rates}.

\begin{remark}\label{rem:comparison}
For comparison, we give a bound obtained directly from
\cite[Theorem 2.3(ii)]{FKS}. Set
\[
\mathcal V_n=\sum_{k=1}^nV_{n,k},\qquad
\mathcal B_n=\left(\sum_{k=1}^nh_{n,k}^2\right)^{1/2}.
\]
Since $V_{n,k}$ is predictable,
\[
\mathcal V_n-\mathsf K_n=\sum_{k=1}^n\Gamma_{n,k},\qquad
\E(\Gamma_{n,k}\mid\F_{k-1})=0.
\]
By~\eqref{eq:Gamma_total_bound}, \eqref{eq:martingale_array} and the
conditional Hoeffding inequality, for $q\ge1$,
\[
\norm{\mathcal V_n-\mathsf K_n}_q\le C_{a,d}\sqrt q\,\mathcal B_n,
\qquad \norm{T_n}_q\le C_{a,d}\sqrt q.
\]
Thus the third-order term in~\cite[equation (2.7)]{FKS} and the
Cauchy--Schwarz inequality give
\begin{align*}
\Wp\bigl(\Lw(T_n),N(0,\mathsf K_n)\bigr)
&\le C_{a,d}p^{2/3}
\left(\sum_{k=1}^nw_{n,k}^{3/2}+\mathcal B_n\right)^{1/3}
\le C_{a,d}p\mathcal B_n^{1/3},\\
\sum_{k=1}^nw_{n,k}^{3/2}
&\le\left(\sum_{k=1}^nw_{n,k}\right)^{1/2}
\left(\sum_{k=1}^nw_{n,k}^2\right)^{1/2}\le\mathcal B_n.
\end{align*}
By~\eqref{eq:weight_asymptotics} and~\eqref{eq:hnk_def},
\[
\mathcal B_n\lesssim
\begin{cases}
n^{-1/2},&a<1/4,\\
n^{-1/2}(\log n)^{1/2},&a=1/4,\\
n^{-(1-2a)},&1/4<a<1/2,\\
(\log n)^{-1},&a=1/2.
\end{cases}
\]
For $a<0$, we also use $h_{n,k}\lesssim n^{-1}$.
Together with~\eqref{eq:gaussian_cov_compare}, this yields
\[
\Wp\bigl(\Lw(T_n),N(0,I_d)\bigr)\le C_{a,d}p
\begin{cases}
n^{-1/6},&a<1/4,\\
n^{-1/6}(\log n)^{1/6},&a=1/4,\\
n^{-(1-2a)/3},&1/4<a<1/2,\\
(\log n)^{-1/3},&a=1/2.
\end{cases}
\]
These are upper bounds from the stated application of the general theorem.
For $a=0$, the increments are i.i.d. and $\mathcal V_n=I_d$;
\cite[Corollary 2.1]{FKS} directly gives $Cp\log(en)/\sqrt n$.
\end{remark}

\end{document}